\documentclass{amsart}

\usepackage{amsfonts}

\newtheorem{theorem}{Theorem}[section]
\newtheorem{lemma}[theorem]{Lemma}
\theoremstyle{definition}
\newtheorem{definition}[theorem]{Definition}
\newtheorem{example}[theorem]{Example}

\theoremstyle{remark}
\newtheorem{remark}[theorem]{Remark}
\numberwithin{equation}{section}
\theoremstyle{plain}

\newtheorem{corollary}[theorem]{Corollary}

\newtheorem{proposition}[theorem]{Proposition}

\newtheorem{observation}[theorem]{Observation}
\copyrightinfo{2001}{enter name of copyright holder}

\newcommand{\N}{\mathbb{N}}

\begin{document}
\title[Isomorphisms of noncommutative domain algebras II]{Isomorphisms of
noncommutative domain algebras II}
\author{Alvaro Arias}
\address{Department of Mathematics\\
University of Denver}
\email{aarias@math.du.edu}
\urladdr{http://www.math.du.edu/\symbol{126}aarias}
\author{Fr\'{e}d\'{e}ric Latr\'{e}moli\`{e}re}
\address{Department of Mathematics\\
University of Denver}
\email{frederic@math.du.edu}
\urladdr{http://www.math.du.edu/\symbol{126}frederic}
\date{October 2010}

\begin{abstract}
This paper extends the results of the previous work of the authors on the classification on noncommutative domain algebras up to completely isometric isomorphism. Using Sunada's classification of Reinhardt domains in $\mathbb{C}^n$, we show that aspherical noncommutative domain algebras are isomorphic if and only if their defining symbols are equivalent, in the sense that one can be obtained from the other via permutation and scaling of the free variables. Our result also shows that the automorphism groups of aspherical noncommutative domain algebras consists of a subgroup of some finite dimensional unitary group. We conclude by illustrating how our methods can be used to extend to noncommutative domain algebras some results from analysis in $\mathbb{C}^n$ with the example of Cartan's lemma.
\end{abstract}

\maketitle

\section{Introduction}

\bigskip Noncommutative domain algebras, introduced by Popescu \cite{Popescu08}, provide a generalization of noncommutative disk algebras and serve as universal operator algebras for a large class of noncommutative domains, i.e. noncommutative analogues of domains in $\mathbb{C}^n$. This paper extends the results of the previous work of the authors \cite{Arias09} on the classification on noncommutative domain algebras up to completely isometric isomorphism. Our present work uses the classification of Reinhardt domains by Sunada in \cite{Sunada78} to provide a complete classification a large class of noncommutative domain algebras in terms of their defining symbol. This class consists of the aspherical noncommutative domain algebras whose symbol is polynomial, as we shall define in the first section of this paper. 

\bigskip Let us introduce the objects of interest in this paper, as well as the notations we will use throughout our exposition. Popescu noncommutative domains are defined by means of a symbol, which is a special type of formal power series. Let $\mathbb{F}_n^+$ be the free semigroup on $n$ generators $g_1,\ldots,g_n$ and identity $g_0$. If $X_1,\ldots,X_n$ belong to any ring, and if $\alpha\in\mathbb{F}_n^+$ is written as $\alpha=g_{i_1}\cdots g_{i_n}$, then we write $X_\alpha = X_{i_1}\cdots X_{i_n}$. With these notations, a free formal power series $f=\sum_{\alpha\in\mathbb{F}_n^+} a_\alpha X_\alpha$ with real coefficients $a_\alpha$ ($\alpha\in\mathbb{F}_n^+$) is \emph{regular positive} if:
\begin{equation}
\left\{ 
\begin{array}{l}
a_{g_{0}}=0\text{,} \\ 
a_{g_{i}}>0\text{ if }i=1,\ldots ,n\text{,} \\ 
a_\alpha\geq 0 \;\textrm{for all $\alpha\in\mathbb{F}_n^+$}\\
\sup_{n\in \mathbb{N}^{\ast }}\left( \left\vert \sum_{\left\vert \alpha
\right\vert =n}a_{\alpha }^{2}\right\vert ^{\frac{1}{n}}\right) <\infty 
\text{.} 
\end{array} 
\right.   \label{PositiveRegularFormal}
\end{equation} 

\bigskip Let us be given a regular positive free formal power series $f=\sum_{\alpha\in\mathbb{F}_n^+} a_\alpha X_\alpha$ in $n$ indeterminates, and a Hilbert space $\mathcal{H}$. We define the noncommutative domain:
$$\mathcal{D}_f(\mathcal{H}) = \{ (T_1,\ldots,T_n) \in \mathcal{B}(\mathcal{H}) : \sum_{\alpha\in\mathbb{F}_n^+} a_\alpha T_\alpha T^*_\alpha \leq 1 \}$$
 where $\mathcal{B}(\mathcal{H})$ is the Von Neumann algebra of all bounded linear operators on $\mathcal{H}$ and $1$ is its identity. The interior of $\mathcal{D}_f(\mathbb{C}^k)$ will be denoted by $\mathbb{D}_f^k$.
Popescu proved in \cite{Popescu08} that certain weighted shifts on the full Fock space $\ell^2(\mathbb{F}_n^+)$ of $\mathbb{C}^n$ provide a model for all $n$-tuples in all domains $\mathcal{D}_f(\mathcal{H})$.  Specifically, as shown in \cite{Popescu08}, one can
find positive real numbers $\left( b_{\alpha }\right) _{\alpha \in \mathbb{F} 
_{n}^{+}}$ and define linear operators $W_{i}^f$ on $\ell ^{2}\left( \mathbb{F} 
_{n}^{+}\right) $ such that: 
\begin{equation*}
W_{i}^f\delta _{\alpha }=\sqrt[2]{\frac{b_{\alpha }}{b_{g_{i}\alpha }}}\delta
_{g_{i}\alpha }\text{,}
\end{equation*} 
where $\{\delta _{\alpha }:\alpha\in\mathbb{F}_n^+\}$ is the canonical basis of $\ell^2(\mathbb{F}_n^+)$. Then $(W_1^f,\ldots,W_n^f)\in\mathcal{D}_f(\ell^2(\mathbb{F}^+_n))$ and the coefficients $(a_\alpha)_{\alpha\in\mathbb{F}_n^+}$ can be recovered from the coefficients $(b_\alpha)_{\alpha\in\mathbb{F}_n^+}$. We then define:
$$\mathcal{A}(\mathcal{D}_f) = \overline{ \mathrm{span} \{ W_\alpha  : \alpha\in\mathbb{F}_n^+ \} }$$
i.e. $\mathcal{A}(\mathcal{D}_f)$ is the norm closure in $\mathcal{B}(\mathcal{H})$ of the algebra generated by $W_1^f,\ldots,W_n^f$. The fundamental property of this algebra is that:
\begin{proposition}[Popescu]
Let $f=\sum_{\alpha\in\mathbb{F}_n^+} a_\alpha X_\alpha$ be a regular positive free power series in $n$ indeterminates. Let $\mathcal{H}$ be a Hilbert space. Let $(T_1,\ldots,T_n)\in\mathcal{D}_f(\mathcal{H})$. Then there exists a (necessarily unique) completely contractive unital algebra morphism $\Phi:\mathcal{A}(\mathcal{D}_F)\longrightarrow \mathcal{B}(\mathcal{H})$ such that $\Phi(W_j^f)=T_j$ for $j=1,\ldots,n$.
\end{proposition}

The algebra $\mathcal{A}(\mathcal{D}_f)$ is the \emph{noncommutative domain algebra of symbol $f$}. When $f=X_{1}+\cdots +X_{n}$, the algebra $\mathcal{A}\left( \mathcal{D}_{f}\right) $ is the disk algebra in $n$ -generators. In \cite{Arias09}, we used techniques of complex analysis on domains in $\mathbb{C}^n$ to study the isomorphism problem for noncommutative domain algebras. In our context, the category \textrm{NCD} of noncommutative domain algebras consists of the algebras $\mathcal{A}(\mathcal{D}_f)$ for all positive, regular $n$-free formal power series $f$ for objects, and \emph{completely isometric} unital algebra isomorphisms for arrows. We then constructed  for each $k \in \mathbb{N}$ a contravariant functor $\mathbb{D}^k$ from {NCD} to the category $\mathrm{HD_{kn^2}}$ of connected open subsets of $\mathbb{C}^{kn^2}$ with holomorphic maps. Given a regular positive formal $n$-free formal power series $f$, the functor $\mathbb{D}^k$ associates to $\mathcal{A}(\mathcal{D}_f)$ the domain $\mathbb{D}_f^k$. The key observation of \cite{Arias09} is that any isomorphism $\Phi:\mathcal{A}(\mathcal{D}_f)\rightarrow\ \mathcal{A}(\mathcal{D}_g)$ in \textrm{NCD} gives rise to a biholomophic map $\widehat{\Phi_k} = \mathbb{D}^k(\Phi)$ from $\mathbb{D}_g^k$ onto $\mathbb{D}_f^k$, and this construction is functorial (contravariant).
Thus, $\mathbb{D}^k$ give fundamental invariants of noncommutative domain algebras. We showed in \cite[Theorem 3.18]{Arias09} that if  $\Phi:\mathcal{A}(\mathcal{D}_f)\rightarrow\mathcal{A}(\mathcal{D}_g)$ is an isomorphism in \textrm{NCD} such that $\widehat{\Phi_1}(0)=0$ then there is an invertible matrix $M\in M_{n\times n}(\mathbb{C})$ such that $\left[ \begin{array}{c}W_1^g\\ \vdots\\ W_n^g \end{array} \right] = (M\otimes 1_{\ell(\mathbb{F}_n^+)}) \left[ \begin{array}{c}W_1^f\\ \vdots\\ W_n^f \end{array} \right]$. 
We could then use this result to characterize the disk algebra among all noncommutative domain algebras \cite[sec 4.3]{Arias09}: $\mathcal{A}(\mathcal{D}_f)$ is isomorphic to the $n$-disk algebra in \textrm{NCD}  if and only if $f=X_1+\ldots+X_n$ after possible rescaling (i.e. replacing each indeterminate with a multiple of itself). We could also distinguish between the instructive examples $f=X_1+X_2+X_1X_2$ and $g=X_1+X_2+\frac{1}{2}X_1X_2+\frac{1}{2}X_2X_1$ in \textrm{NCD} \cite[sec. 4.2]{Arias09}.

\bigskip In this paper, we extend our results to characterize a large class of noncommutative domain algebras in terms of their symbol. For instance, we strengthen \cite[sec. 4.2]{Arias09} by showing that, given  $f=X_1+X_2+X_1X_2$, the algebra $\mathcal{A}(\mathcal{D}_f)$ is isomorphic in \textrm{NCD} to $\mathcal{A}(\mathcal{D}_g)$ if and only if, after permuting and rescaling the indeterminates in $g$, we have $f=g$. 

\bigskip In \cite[Theorem 4.5]{Popescu11}, Popescu proved that two noncommutative domain algebras are isomorphic in \textrm{NCD} if and only if their corresponding noncommutative domains are biholomorphic, in the generalized sense of Popescu. Thus, our result extends to the biholomorphic classification of noncommutative domains in terms of their defining symbols.
We note that Popescu's result applies to a larger class of domain algebras and noncommutative domains associated with the Berezin transform. We refer to \cite[Theorem 4.5]{Popescu11} for details.

\bigskip Thus, this paper continues the program initiated in \cite{Arias09} with three new results. First, we show that Sunada's classification of Reinhardt domains \cite{Sunada78} allows one to prove that for a large class of noncommutative domains, the only possible automorphisms are linear. Second, we show that
for that class of domains, the isomorphism problem is fully solved, as
isomorphisms correspond to a simple form of equivalence on regular positive $n$-free formal power series . As a third matter, we show that the same
techniques allow us to obtain a generalization of the Cartan's lemma to noncommutative domain by different means than \cite{Popescu10,Popescu11}. We expect that many
results of complex analysis in several variables can be generalized to
noncommutative domain algebras in a similar manner.

\section{aspherical noncommutative domain algebras}

This first section applies Sunada's classification of Reinhardt domains to
the classification of noncommutative domain algebras. In \cite{Arias09}, we showed how to
classify all noncommutative domain algebras isomorphic to the disk algebras.
In this paper, we will provide a complete classification for the class of polynomial aspherical domains as defined below.

\begin{definition}
Let $f,g$ be two regular positive free power series in $n$ indeterminates. We say that $f$ and $g$ are \emph{permutation-rescaling equivalent} if there exists a permutation $\sigma$ of $\{1,\ldots,n\}$ and a nonzero positive numbers $\lambda_1,\ldots,\lambda_n$ such that $g(\lambda_1 X_{\sigma(1)},\ldots,\lambda_n X_{\sigma(n)})=f(X_1,\ldots,X_n)$ where $X_1,\ldots,X_n$ are the indeterminates.
\end{definition}

One checks easily that permutation-rescaling equivalence is an equivalence relation of regular positive free power series. It was shown in \cite[Lemma 4.4]{Arias09} that if $f$ and $g$ are regular positive free formal power series in $n$ indeterminates, and $f$ and $g$ are permutation-rescaling equivalent, then $\mathcal{A}(\mathcal{D}_f)$ and $\mathcal{A}(\mathcal{D}_g)$ are isomorphic.

\begin{definition}
Let $f$ be a positive regular $n$-free formal
power series. Then $f$ is spherical if it is permutation-rescaling to some $g$ such that: 
\begin{equation*}
\mathbb{D}_g^1 = \{ (z_1,\ldots,z_n) \in \mathbb{C}^n : |z_1|^2+\cdots+|z_n|^2 < 1 \} \textrm{.}
\end{equation*} 
\end{definition}

\begin{definition}
Let $f$ be a positive regular $n$-free formal power series. Then $f$ is aspherical if $f$ is not spherical.
\end{definition}

\begin{definition}
 A free power series $f$ is a free polynomial if, writing $f=\sum_{\alpha\in\mathbb{F}^+_n} a_\alpha X_\alpha$, the set $\{\alpha\in\mathbb{F}_n^+:a_\alpha\not=0\}$ is finite. We shall write $n$-free polynomial for a free polynomial in $n$-indeterminates. 
\end{definition}

\begin{definition}
The noncommutative domain $\mathcal{A}\left( \mathcal{D}_{f}\right) $ is
aspherical when $f$ is aspherical, and polynomial when $f$ is a free polynomial.
\end{definition}

\bigskip It should be noted that there are domains which are not aspherical but
not isomorphic to the disk algebras. Let us give a few examples to illustrate our definition.

\begin{example}
Let $f=\frac{1}{2}X_1+\frac{1}{2}X_2+\frac{1}{2}X_1^2+\frac{1}{2}X_2^2+X_1X_2$. Then $\mathbb{D}_f^1$ is the open unit ball of $\mathbb{C}^2$. Yet, by \cite[Theorem 4.7]{Arias09}, we know that $\mathcal{A}(\mathcal{D}_f)$ is not a disk algebra.
\end{example}

\begin{example}
Let $f=\frac{1}{2}X_1+\frac{1}{2}X_2+\frac{1}{2}(X_1+X_2)^2$. Using \cite[Theorem 4.7]{Arias09}, we see again that $\mathcal{A}(\mathcal{D}_f)$ is not isomorphic to a disk algebra in \textrm{NCD}. However, there is a completely bounded isomorphism from the disk algebra onto $\mathcal{A}(\mathcal{D}_f)$. This illustrates the importance of our choice of isomorphisms as \emph{completely isometric} unital algebra isomorphisms.
\end{example}

\begin{example}\label{silly}
Let $f=X_1+X_2+3X_1 X_2$, $g = 2X_1+X_2+6X_2 X_1$ and $h=X_1+2X_2+X_1^2$. All three symbols are polynomial and aspherical. We will show in Theorem (\ref{Main}) that $\mathcal{A}(\mathcal{D}_f)$ and  $\mathcal{A}(\mathcal{D}_g)$ are isomorphic, but  $\mathcal{A}(\mathcal{D}_f)$ and  $\mathcal{A}(\mathcal{D}_g)$ are not.
\end{example}

It is of course quite easy to produce examples of polynomial aspherical symbols, so our work will apply to a large class of examples.

\bigskip It is immediate, by definition, that:

\begin{observation}
Let $f$ be a positive regular free formal power series in $n$ indeterminates. Then $ 
\mathbb{D}_f^1 $ is a bounded Reinhardt domain in 
$\mathbb{C}^{n^2}$.
\end{observation}

\bigskip Now, Sunada in \cite[Theorem A]{Sunada78} shows that up to rescaling-permutation of the canonical basis vectors of $\mathbb{C}^n$, all bounded Reinhardt domains can be written in a normalized form. We cite this theorem with the necessary notations added in the statement.

\begin{theorem}[Sunada, Theorem A]\label{TheoremA}
Let $D$ be a bounded Reinhardt domain in $\mathbb{C}^n$. Up to applying to $D$ a map of the form $(z_1,\ldots,z_n) \in \mathbb{C}^n \mapsto (\lambda_1 z_{\sigma(1)},\ldots,\lambda_n z_{\sigma(n)})$ with $\sigma$ some permutation of $\{1,\ldots,n\}$ and $\lambda_1,\ldots,\lambda_n$ some nonzero complex numbers, $D$ can be described as follows. There exists integers $0 = n_0 < n_1 < \ldots < n_s =n$, integers $p,r\in\{1,\ldots,n\}$ with $n_r=p$ and a bounded Reinhard domain $D_1$ in $\mathbb{C}^{n-p}$ such that, if for any $z=(z_1,\ldots,z_n)\in\mathbb{C}^n$ we set $\mathbf{z}_j=(z_{n_{j-1}+1},\ldots,z_{n_j})$ and $r=n-p$, then:
\begin{itemize}
\item $D \cap (\mathbb{C}^p \times \{0\})= \{ (\mathbf{z}_1,\ldots,\mathbf{z}_s) \in \mathbb{C}^n : |\mathbf{z}_1| < 1 \wedge \cdots \wedge |\mathbf{z}_r| < 1 \wedge (\mathbf{z}_{r+1},\ldots,\mathbf{z}_{s})=(0,\ldots,0)\}$
\item $\{0\}\times D_1 = D \cap (\{0\} \times \mathbb{C}^{n-p})$,
\item We have:
\begin{eqnarray*}
D &=& \left\{ (\mathbf{z}_1,\ldots,\mathbf{z}_s) : |\mathbf{z}_1|<1\wedge \ldots \wedge |\mathbf{z}_r|<1 \wedge \right.\\
& & { } \left. \left(\frac{\mathbf{z}_{r+1}}{\prod_{j=1}^r (1-|\mathbf{z}_j|^2)^{q_{r+1,j}}},\ldots,\frac{\mathbf{z}_{s}}{\prod_{j=1}^r (1-|\mathbf{z}_j|^2)^{q_{s,j}}}\right) \in D_1 \right\}
\end{eqnarray*}
where $q_{r+1,1},\ldots,q_{s,r}$ are positive nonzero real numbers distinct from $1$.
\end{itemize}
\end{theorem}

We also observe that if $f$ is a regular positive $n$-free power series, then $\mathbb{D}_f^1 \cap (\mathbb{C}^p \times \{0\})=\mathbb{D}_g^1$ for the power series $g$ obtained by evaluating $X_{p+1},\ldots,X_n$ to $0$. One readily checks that $g$ is a regular positive $p$-free power series. This simple observation will prove useful.

Our strategy to prove our main theorem of this section, Theorem (\ref{main}), is as follows. Let $f$ be a regular positive aspherical $n$-free polynomial.

\begin{enumerate}

\item Since $f$ is aspherical, $n_1<n$ and thus $r<s$ in Theorem (\ref{TheoremA}) for $\mathbb{D}_f^1$.

\item We first show that $\mathbb{D}_f^1$ can not be a product of bounded Reinhardt domains. This implies that $r \in \{0,1\}$ in Theorem (\ref{TheoremA}) for $\mathbb{D}_f^1$.
\item We then show that if $n=2$, i.e. $f$ has only two indeterminates, then $\mathbb{D}_f^1$ can not be a Thullen domain \cite{Thullen31}. This implies that $r=0$ in Theorem (\ref{TheoremA}) for $\mathbb{D}_f^1$ for an arbitrary $n$.

\item We then conclude, using (\cite[Corollary 2, sec. 6. p. 126]{Sunada78}), that the automorphism group of $\mathbb{D}_f^1$ fixes the origin.
\end{enumerate}

To implement this approach, we start with a few lemmas.

\begin{lemma}\label{noproduct}
Let $f$ be a positive regular $n$-free formal polynomial. Then $\mathcal{D} 
_{f}\left( \mathbb{C}\right) $ is not a Cartesian product of bounded domains.
\end{lemma}

\begin{proof}
Up to rescaling, we assume that $a_\alpha = 1$ for all words $\alpha$ of length 1. Assume that there are two domains $D$ and $D^{\prime }$, respectively in $ 
\mathbb{C}^{p}$ and $\mathbb{C}^{q}$ with $p+q=n$, $p,q>0$ and such that: 
\begin{equation*}
\mathbb{D}_{f}^1 =D\times D^{\prime }\text{.}
\end{equation*} 
Let $\left( z_{1},\ldots ,z_{p}\right) \in D$ be a boundary point. Then $ 
\left( z_{1},\ldots ,z_{p},0,\ldots ,0\right) \in \mathbb{C}^{n}$ is also a
boundary point of $\mathbb{D}_f^1 $. Hence: 
\begin{equation*}
\sum_{\alpha \in \mathbb{F}_{p}^{+}}a_{\alpha }\left\vert z_{\alpha
}\right\vert ^{2}=1\text{.}
\end{equation*} 
Let $\left( w_{1},\ldots ,w_{q}\right) \in D^{\prime }\backslash \left\{
0\right\} $. Then $\left( z_{1},\ldots ,z_{p},w_{1},\ldots ,w_{q}\right) \in 
\mathbb{D}_f^1 $, so (if we let $z_{p+1}=w_{1},\ldots,z_n = w_{q}$): 
\begin{eqnarray*}
1 &\geq &\sum_{\alpha \in \mathbb{F}_{n}^{+}}a_{\alpha }\left\vert z_{\alpha
}\right\vert ^{2}\geq \sum_{\alpha \in \mathbb{F}_{p}^{+}}a_{\alpha
}\left\vert z_{\alpha }\right\vert ^{2}+\left\vert w_{1}\right\vert ^{2}+\cdots +\left\vert w_{q}\right\vert ^{2} \\
&\geq &1+\left\vert w_{1}\right\vert ^{2}+\cdots +\left\vert
w_{q}\right\vert ^{2}>1
\end{eqnarray*} 
which is a contradiction. Hence $\mathbb{D}^1_f $
is not a product of domains in proper subspaces of $\mathbb{C}^n$.
\end{proof}

Thullen \cite{Thullen31} proved that all bounded Reinhardt domains containing the origin in $\mathbb{C}^2$ are biholomorphic to either the unit ball, the polydisk, a domain whose automorphism group fixes $0$, or to a Thullen domain:
$$\{ (z,w) \in\mathbb{C}^2 : |z|^2+|w|^{2q} \leq 1 \}$$
for some $q \in (0,1)\cup(1,\infty)$. This result, of course, is generalized in Theorem (\ref{TheoremA}). We will start by showing that in two dimensions, Thullen domains are not in the range of the object map of the functor s $\mathbb{D}^1$.

\begin{lemma}\label{nothullen2d}
Let $f$ be an aspherical positive regular $2$-free formal polynomial. Then $ 
\mathcal{D}_{f}\left( \mathbb{C}\right) $ is not a Thullen domain \cite{Thullen31}.
\end{lemma}

\begin{proof}

Assume $\mathcal{D}_{f}\left( \mathbb{C}\right) $ is a Thullen domain \cite{Thullen31}. Up to applying a rescaling and permutation of the free variables of $f$ (\cite[Lemma 4.4]{Arias09}), there exists $q\in \left( 0,1\right) \cup \left( 1,\infty \right) $
such that, writing $\tau :\left( z_{1},z_{2}\right) \mapsto \left\vert
z_{1}\right\vert ^{2}+\left\vert z_{2}\right\vert ^{2 q}-1$, we have:

$$\mathbb{D}_f^1 =\left\{ \left( z,w\right) \in 
\mathbb{C}^{2}:\tau \left( z,w\right) < 0\right\} $$ and the boundary $\partial \mathbb{D}_f^1$ of $\mathbb{D}_f^1$ is then $$\left\{ \left(
z_{1},z_{2}\right) \in \mathbb{C}^{2}:\tau \left( z_{1},z_{2}\right)
=0\right\} \textrm{.}$$ 

We shall adopt the standard notation that given $(z_1,z_2)\in\mathbb{C}^2$, we have $z_j=x_j+ i y_j$ for $x_j,y_j\in \mathbb{R}$. We identity $\mathbb{C}^2$ with $\mathbb{R}^4$. Then:
$$\tau: (x_1,y_1,x_2,y_2) \mapsto x_1^2 + y_1^2 + (x_2^2+y_2^2)^q - 1\textrm{.}$$

Let $\nabla$ denote the gradient operator on $\mathbb{R}^{4}$ (i.e., with usual abuse of notations, $\nabla=\left[ \begin{array}{c} \frac{\partial}{\partial x_1} \\ \frac{\partial}{\partial y_1} \\ \frac{\partial}{\partial x_2} \\ \frac{\partial}{\partial y_2} \end{array}\right]$). Then for all $(z_1,z_2) \in \partial \mathbb{D}_f^1)$:

\[
\nabla \tau(z_1,z_2) = \left[ \begin{array}{c} 2x_1 \\ 2y_1 \\ 2q x_2 (x_2^2+y_2^2)^{q-1} \\ 2q y_2 (x_2^2+y_2^2)^{q-1} \end{array} \right] \textrm{.} 
\]

On the other hand, by definition, the boundary of $\mathbb{D}_f^1$ is: $$\left\{
\left( z_{1},z_{2}\right) \in \mathbb{C}^{2}:\sum a_{\alpha }\left\vert
z_{\alpha }\right\vert ^{2}-1=0\right\} \textrm{.}$$ To simplify notations, we shall
introduce the coefficients $\left( c_{n,m}\right) _{n,m\in \N}$
as follows: 
\begin{equation*}
c_{n,m}=\sum \left\{ a_{\alpha }:\alpha \in \mathbb{F}_{2}^{+},\text{ 
\thinspace }\left\vert \alpha \right\vert _{1}=n\wedge \left\vert \alpha
\right\vert _{2}=m\right\}
\end{equation*} 
where $\left\vert \alpha \right\vert _{i}$ is the number of times the
generator $g_{i}$ appears in the word $\alpha $ ($i=1,2$). Thus we can write
for all $\left( z_{1},z_{2}\right) \in \mathbb{C}^{2}$: 
\begin{equation*}
\sum_{\alpha \in \mathbb{F}_{2}^{+}}a_{\alpha }\left\vert z_{\alpha
}\right\vert ^{2}=\sum_{n,m\in \N}c_{n,m}\left\vert
z_{1}\right\vert ^{2n}\left\vert z_{2}\right\vert ^{2m}\text{.}
\end{equation*}

Now, let: 
\begin{equation*}
\rho :\left( z_{1},z_{2}\right) \in \mathbb{C}^{2}\mapsto \sum_{n,m\in 
\N}c_{n,m}\left\vert z_{1}\right\vert ^{2n}\left\vert
z_{2}\right\vert ^{2m}-1\text{.}
\end{equation*} 

We thus have:
$$\rho(x_1,y_1,x_2,y_2) = \sum_{n,m\in\N} c_{n,m} (x_1^2+y_1^2)^n (x_2^2 + y_2^2)^m - 1\textrm{.}$$

Then for all $(z_1,z_2) \in \partial \mathbb{D}_f^1$: 
\begin{eqnarray*}
\nabla \rho(z_1,z_2)&=& \left[ \begin{array}{c} 
\sum_{n,m\in\N,(n,m)\not=(0,0)} c_{n,m} 2nx_1(x_1^2+y_1^2)^{n-1} (x_2^2 + y_2^2)^m \\
\sum_{n,m\in\N,(n,m)\not=(0,0)} c_{n,m} 2ny_1(x_1^2+y_1^2)^{n-1} (x_2^2 + y_2^2)^m \\
\sum_{n,m\in\N,(n,m)\not=(0,0)} c_{n,m} 2mx_2(x_1^2+y_1^2)^{n} (x_2^2 + y_2^2)^{m-1} \\
\sum_{n,m\in\N,(n,m)\not=(0,0)} c_{n,m} 2my_2(x_1^2+y_1^2)^{n} (x_2^2 + y_2^2)^{m-1}
\end{array} \right]\\
&=&
\left[ \begin{array}{c} 
2x_1 \left(c_{1,0}+ \sum_{n,m\in\N,n>1} nc_{n,m} (x_1^2+y_1^2)^{n-1} (x_2^2 + y_2^2)^m \right)\\
2y_1 \left(c_{1,0} + \sum_{n,m\in\N,n>1}nc_{n,m} (x_1^2+y_1^2)^{n-1} (x_2^2 + y_2^2)^m \right)\\
2x_2 \left(c_{0,1} + \sum_{n,m\in\N,m>1}mc_{n,m} (x_1^2+y_1^2)^{n} (x_2^2 + y_2^2)^{m-1} \right)\\
2y_2 \left(c_{0,1} + \sum_{n,m\in\N,m>1}mc_{n,m} (x_1^2+y_1^2)^{n} (x_2^2 + y_2^2)^{m-1} \right)
\end{array} \right]\\
\end{eqnarray*}
The tangent plane in $\mathbb{R}^4$ to a boundary point $\left( z_{1},z_{2}\right) 
$ (where $z_2\not=0$ so that we work at a regular point for $\tau$) of $\mathbb{D}_f^1$ is the orthogonal space to
any normal vector to the boundary of $\mathbb{D}_f^1$ at that point, namely it is the orthogonal of $\nabla\rho \left(
z_{1},z_{2}\right) $, as well as the orthogonal space to $\nabla\tau \left( z_{1},z_{2}\right) $ (see, for instance, \cite[chapter 3]{Krantz}). Thus, these two vectors must be co-linear. In particular, let us focus on the first and third coordinates.
It is therefore necessary that if $x_1\not=0, x_2\not=0$: 
\begin{eqnarray*}
&{}& \left( c_{1,0}+ \sum_{n,m\in\N,n>1} nc_{n,m} (x_1^2+y_1^2)^{n-1} (x_2^2 + y_2^2)^m   \right)  q(x_2^2+y_2^2)^{q-1} \\
&=& c_{0,1} + \sum_{n,m\in\N,m>1}mc_{n,m} (x_1^2+y_1^2)^{n} (x_2^2 + y_2^2)^{m-1}
\text{.}
\end{eqnarray*} 
which is equivalent to:
\begin{eqnarray*}
&{}& \left(c_{1,0} + \sum_{n,m\in\N,n>1} nc_{n,m} \vert z_1 \vert ^{2n-2} \vert z_2 \vert ^{2m} \right) q \vert z_2 \vert ^{2q-2} \\
&=& c_{0,1} + \sum_{n,m\in\N,m>1}mc_{n,m} \vert z_1 \vert^{2n} \vert z_2 \vert ^{2m-2} \textrm{.} 
\end{eqnarray*}
Now, since $\left( z_{1},z_{2}\right) $ is on the boundary of the Thullen
domain $\mathbb{D}_f^1$, we have $\tau \left(
z_{1},z_{2}\right) =0$ i.e. $\left\vert z_{1}\right\vert ^{2}=1-\left\vert
z_{2}\right\vert ^{2q}$. Hence: 
\begin{eqnarray*}
&{}&\left( c_{1,0}+\sum_{n,m\in\N, n>1}  nc_{n,m}\left( 1-\left\vert z_{2}\right\vert ^{2q}\right)
^{n-1}\left\vert z_{2}\right\vert ^{2m}\right) q\left\vert z_{2}\right\vert
^{2q-2}\\
&=&c_{0,1}+\sum_{n,m\in \N, m>1} mc_{n,m}\left( 1-\left\vert z_{2}\right\vert ^{2q}\right)
^{n}\left\vert z_{2}\right\vert ^{2m-2}\text{.}
\end{eqnarray*} 
This identity must be valid for all $\left( z_{1},z_{2}\right) $ on the
boundary of $\mathbb{D}_f^1$ except when $z_{2}=0$. Thus it is true on a neighborhood of $0$, except at $0$. Hence both sides of this identity must be continuous at $0$ since the right hand side is as a polynomial in $\mathbb{R}^2$. This precludes that $q<1$. Hence $q>1$ and we get at the limit when $z_2\rightarrow 0$ that $0=c_{0,1}$, which contradicts the definition of $f$ regular positive. So $\mathcal{D}_{f}\left( \mathbb{C}\right) $ is not a Thullen
domain.
\end{proof}

\begin{remark}\label{thullen}

An important observation is that Thullen domains are a special case of the domains described in Theorem (\ref{TheoremA}). Indeed, assume, in the notations of \ref{TheoremA}, that $n=2$, $s=2$, $r=1$. Then $D_1$, which is a bounded Reinhard domain in $\mathbb{C}$ which we can put in standard form, is just the unit disk in $\mathbb{C}$, so:
\begin{eqnarray}
D &=& \{ (z_1,z_2)\in\mathbb{C}^2 : |z_1|<1\wedge |z_2|<(1-|z_1|^2)^q \} \\
&=& \{ (z_1,z_2)\in\mathbb{C}^2 : |z_1|^2+|z_2|^{2/q} < 1 \}
\end{eqnarray}

which is to say that $D$ is a Thullen domain \cite{Thullen31}. In particular, if $n$ is now arbitrary, $s>r>1$, then the intersection of $D$ with the plane spanned by the first and $(p+1)^{th}$ canonical basis vectors of $\mathbb{C}^n$ is a Thullen domain.

\end{remark}

\begin{theorem}\label{main}
Let $f$ be a aspherical regular positive $n$-free polynomial. Then the
automorphism group of $\mathcal{D}_{f}\left( \mathbb{C}\right) $ fixes $0$.
\end{theorem}

\begin{proof}

By \cite[Proposition 3.11]{Arias09}, $\mathbb{D}_f^1$ is a bounded Reinhardt domain in $\mathbb{C}^n$.

\bigskip Up to replacing $f$ by a permutation-scaling equivalent symbol, we can assume that $\mathbb{D}_f^1$ is in standard form, i.e. of the form given in Theorem (\ref{TheoremA}). We recall from \cite[Lemma 4.4]{Arias09} that if two symbols are permutation-scaling equivalent, then their associated noncommutative domain algebras are isomorphic.

\bigskip We shall now use the notations of Theorem (\ref{TheoremA}) applied to $\mathbb{D}_f^1$.

\bigskip By definition of aspherical, $\mathbb{D}_f^1$ is not the open unit ball of $\mathbb{C}^n$, so $n_1<n$, $r<s$ and $p<n$.

\bigskip Assume now that $r\geq 1$. Then $\mathbb{D}_f^1\cap (\mathbb{C}^{p}\times\{0\})$ is a product of open unit balls. Yet, if $g$ is obtained from $f$ by mapping $X_{p+1},\ldots,X_n$ by $0$, then $g$ is a regular positive $p$-free polynomial such that $\mathbb{D}_g^1 \times {0} = \mathbb{D}_f^1 \cap (\mathbb{C}^p \times\{0\})$ by construction. By Lemma (\ref{noproduct}), $\mathbb{D}_g^1$ is not a proper product, so it must be at most one unit ball. Hence, $r=1$.

\bigskip Therefore, in general, $r\in \{0,1\}$. Assume that $r=1$, so $1\leq p<n$ (since $r<s$ as well). Let $h$ be the regular formal $2$-free polynomial obtained from $f$ by mapping $X_2,\ldots,X_p,X_{p+2},\ldots,X_n$ to $0$ (once again it is immediate to check that indeed $h$ is a regular positive polynomial in indeterminates $X_1,X_{p+1}$). In particular, observe that $\mathbb{D}_h^1$ is the intersection of $\mathbb{D}_f^1$ by the plane spanned by the first and $(p+1)^{th}$ canonical basis vector in $\mathbb{C}^n$. By Remark (\ref{thullen}), $\mathbb{D}_h^1$ is a Thullen domain. By Lemma (\ref{nothullen2d}), this is impossible. So we have reached a contradiction.

\bigskip Hence $r=0$. By \cite[Corollary 2 of Theorem B, sec 6. p.126]{Sunada78}, all automorphisms of $\mathbb{D}_f^1$ fix $0$. Our theorem is proven.

\end{proof}

\bigskip Applying \cite[Theorem 3.18]{Arias09}, we get the following important result:

\begin{theorem}\label{Main}
Let $f$ and $g$ be regular positive $n$-free polynomials, with $f$
aspherical. Then if $\Phi :\mathcal{A}\left( \mathcal{D}_{f}\right)
\rightarrow \mathcal{A}\left( \mathcal{D}_{g}\right) $ is an isomorphism,
then $\widehat{\Phi }(0)=0$ and therefore, there exists $M\in M_{n\times
n}\left( \mathbb{C}\right) $ such that: 
\begin{equation*}
\left[ 
\begin{array}{c}
\Phi (W_{1}^{f}) \\ 
\vdots \\ 
\Phi (W_{n}^{f}) 
\end{array} 
\right] =(M\otimes 1_{\ell^2({\mathbb{F}_n})})\left[ 
\begin{array}{c}
W_{1}^{g} \\ 
\vdots \\ 
W_{n}^{g} 
\end{array} 
\right] \text{.}
\end{equation*}
\end{theorem}

\begin{proof}
Let $\lambda =\widehat{\Phi }(0)$. Assume $\lambda \not=0$. Let $j\in
\left\{ 1,\ldots ,n\right\} $ such that $\lambda _{j}\not=0$ where $\lambda
=\left( \lambda _{1},\ldots ,\lambda _{n}\right) $. Let $\rho $ be the
linear transformation of $\mathbb{C}^{n}$ given by the matrix $\left[ 
\begin{array}{ccc}
1_{j-1} &  &  \\ 
& -1 &  \\ 
&  & 1_{n-j} 
\end{array} 
\right] $ in the canonical basis of $\mathbb{C}^{n}$ (where $1_{k}$ is the
identity of order $k$ for any $k\in \mathbb{N}$). Then $\rho (\lambda
)\not=\lambda $. By \cite[Lemma 4.4]{Arias09}, there exists a unique automorphism $\widetilde{\rho }$ of $\mathcal{A}\left( \mathcal{D} 
_{g}\right) $ such that $\mathbb{D}^1(\widetilde{\rho})=\rho$ since $\rho$ only rescale the coordinates. Let $\tau =\Phi ^{-1}\circ \widetilde{\rho }\circ \Phi $: by
construction, $\tau $ is an automorphism of $\mathcal{A}\left( \mathcal{D} 
_{f}\right) $, and moreover $\widehat{\tau }\left( 0\right) =\widehat{\Phi } 
^{-1}\left( \rho \left( \lambda \right) \right) \not=0$ since $\widehat{\Phi 
}^{-1}$ is injective. Thus $\widehat{\tau }$ is an automorphism of $\mathcal{ 
D}_{f}\left( \mathbb{C}\right) $ which does not fix $0$. By Theorem (\ref{main}),
since $f$ is aspherical, this is a contradiction. Hence $\widehat{\Phi } 
(0)=0 $. The theorem follows from \cite[Theorem 3.18]{Arias09}.
\end{proof}

\begin{corollary}
Let $f,g$ be aspherical regular positive $n$-free polynomials. Then there
exist an isomorphism $\Phi $ from $\mathcal{A}\left( \mathcal{D}_{f}\right) $
onto $\mathcal{A}\left( \mathcal{D}_{g}\right) $ if and only if there exists an invertible matrix $M\in M_{n\times n}\left( \mathbb{C}\right) $ such that: 
\begin{equation*}
\left[ 
\begin{array}{c}
W_{1}^{f} \\ 
\vdots \\ 
W_{n}^{f} 
\end{array} 
\right] =(M\otimes 1_{\ell^2(\mathbb{F}_n)})\left[ 
\begin{array}{c}
W_{1}^{g} \\ 
\vdots \\ 
W_{n}^{g} 
\end{array} 
\right] \text{.}
\end{equation*}
\end{corollary}

In particular, the automorphism group of aspherical polynomial noncommutative domain algebras consists only of invertible linear transformation on the generators, in contrast with the disk algebras which have the full automorphism group of the unit ball as a normal subgroup of automorphism \cite{Pitts98}. We also contrast this result to the computation of other nontrivial automorphism groups associated to the disk algebras in \cite{AL10}.

\begin{corollary}
Let $f$ be an aspherical regular positive $n$-free polynomial. The automorphism group of $\mathcal{A}(\mathcal{D}_f)$ is a subgroup of the unitary group $U(n)$.
\end{corollary}

We can now see immediately that Theorem (\ref{Main}) solves Example (\ref{silly}).

\section{Classification of polynomial aspherical noncommutative domains algebras}

\bigskip This section establishes an explicit equivalence on aspherical
regular positive $n$-free polynomials which corresponds to isomorphism of
the associated noncommutative domain algebra. This result generalizes \cite[sec 4]{Arias09}

We defined in \cite{Arias09} dual maps associated to an isomorphism in \textrm{NCD} as follows. Fix $f$ and $g$ two regular, positive $n$-free formal power series, and let $\Phi:\mathcal{A}(\mathcal{D}_f)\rightarrow \mathcal{A}(\mathcal{D}_g)$ be an isomorphism in \textrm{NCD}. Let $k\in\mathbb{N}$, $k>0$. For any $\lambda=(\lambda_1,\cdots,\lambda_n) \in \mathbb{D}^k_f \subseteq (M_{k\times k}(\mathbb{C}))^n$, there exists a unique completely contractive representation $\pi_\lambda : \mathcal{A}(\mathcal{D}_g) \rightarrow M_{n\times n}(\mathbb{C})$ such that $\pi(W_j^g)=\lambda_j$ for $j=1,\ldots,n$. Now, $\rho = \pi_\lambda\circ \Phi$ is a completely contractive representation of $\mathcal{A}(\mathcal{D}_f)$ on $M_{n\times n}(\mathbb{C})$, so there exists $(\mu_1,\ldots,\mu_n)\in\mathbb{D}_f^k$ such that $\rho(W_j^f)=\mu_j$ for $j=1,\ldots,n$. We set $\widehat{\Phi_k}(\lambda_1,\ldots,\lambda_n)=(\mu_1,\ldots,\mu_n)$. We showed in \cite{Arias09} that these maps are biholomorphic and that this construction is functorial and contravariant.

\bigskip We start with a combinatorial result.

\begin{proposition}\label{partition}
Let $f=\sum_{\alpha\in\mathbb{F}_n^+}a_\alpha^f X_\alpha$ and $g=\sum_{\alpha\in\mathbb{F}_n^+} a_\alpha^g X_\alpha$ be regular positive $n$-free polynomials. Up to rescaling, we assume $a_\alpha^f=a_\alpha^g=1$ for all words $\alpha$ of length 1.
Let $\Phi :\mathcal{A} 
\left( \mathcal{D}_{f}\right) \longrightarrow \mathcal{A}\left( \mathcal{D} 
_{g}\right) $ be an isomorphism such that $\widehat{\Phi }(0)=0$. Let $U= 
\left[ u_{ij}\right] _{1\leq i,j\leq n}\in M_{n\times n}\left( \mathbb{C} 
\right) $ be the unitary matrix such that for all $i\in \left\{ 1,\ldots
,n\right\} $: 
\begin{equation}\label{unitarymap}
\Phi \left( W_{i}^{f}\right) =\sum_{j=1}^{n}u_{ij}W_{j}^{g}\text{.}
\end{equation} 
We define the \emph{support function} $s_{\Phi }$ of $\Phi $ by setting, for
any subset $A$ of $\left\{ 1,\ldots ,n\right\} $: 
\begin{equation*}
s_{\Phi }\left( A\right) =\left\{ j\in \left\{ 1,\ldots ,n\right\} :\exists
i\in A\ \ \ u_{ij}\not=0\right\} \text{.}
\end{equation*} 
Then there exists partitions $\left\{ \sigma _{1},\ldots ,\sigma
_{p}\right\} $ and $\left\{ \psi _{1},\ldots ,\psi _{p}\right\} $ of $ 
\left\{ 1,\ldots ,n\right\} $ such that for all $i\in \left\{ 1,\ldots
,p\right\} $ we have $s_{\Phi }\left( \sigma _{i}\right) =\psi _{i}$, $ 
s_{\Phi ^{-1}}\left( \psi _{i}\right) =\sigma _{i}$ and $\left\vert \sigma
_{i}\right\vert =\left\vert \psi _{i}\right\vert $. Moreover, if $A\subseteq
\left\{ 1,\ldots ,n\right\} $ satisfies $s_{\Phi ^{-1}}\circ s_{\Phi }(A)=A$ 
, then: $$A=\bigcup \left\{ \sigma _{i}:i\in \left\{ 1,\ldots ,p\right\}
\wedge \sigma _{i}\subseteq A\right\} \textrm{.}$$
\end{proposition}

\begin{proof}
Since $\Phi $ is an isomorphism such that $\widehat{\Phi }(0)=0$, we
conclude by \cite[Theorem 3.18]{Arias09} that there exists a unitary $U=[u_{ij}]_{1\leq i,j \leq n}\in M_{n\times n}\left( \mathbb{ 
C}\right) $ such that Equality (\ref{unitarymap}) holds for $i\in \left\{ 1,\ldots
,n\right\} $. Note that by definition $s_{\Phi
}(A)=\bigcup\limits_{i\in A}s_{\Phi }\left( \left\{ i\right\} \right) $.
Moreover, for the same reason, there exists $V=\left[ v_{ij}\right] _{1\leq
i,j\leq n}\in M_{n\times n}\left( \mathbb{C}\right) $ such that $\Phi
^{-1}(W_{i}^{g})=\sum_{j=1}^{n}v_{ij}W_{j}^{f}$ and it is immediate that $ 
V=U^{\ast }$ i.e.: 
\begin{equation*}
\Phi ^{-1}\left( W_{i}^{g}\right) =\sum \overline{u_{ji}}W_{j}^{f}\text{.}
\end{equation*} 
This implies that for $A\subseteq \left\{ 1,\ldots ,n\right\} $ we have: 
\begin{eqnarray*}
s_{\Phi ^{-1}}(A) &=&\left\{ j\in \left\{ 1,\ldots ,n\right\} :\exists i\in
A\ \ \ v_{ij}\not=0\right\} \text{ by definition,} \\
&=&\left\{ j\in \left\{ 1,\ldots ,n\right\} :\exists i\in A\ \ \
u_{ji}\not=0\right\} \text{ since }V=U^{\ast}\text{.}
\end{eqnarray*} 
In particular: 
\begin{equation*}
k\in s_{\Phi }\left( \left\{ i\right\} \right) \iff i\in s_{\Phi
^{-1}}\left( \left\{ k\right\} \right) \text{.}
\end{equation*}
Thus $i\in s_{\Phi ^{-1}}\circ s_{\Phi }\left( \left\{ i\right\} \right) $.

We now observe that given any partition $\left\{ \sigma _{1},\ldots ,\sigma
_{n}\right\} $, if we set $\psi _{i}=s_{\Phi }\left( \sigma _{i}\right) $
for all $i\in \left\{ 1,\ldots ,n\right\} ,$ then $\left\vert \psi
_{i}\right\vert \geq \left\vert \sigma _{i}\right\vert $ for all $i\in
\left\{ 1,\ldots ,n\right\} $. Indeed, fix $i\in \left\{ 1,\ldots ,n\right\} 
$. Let $j\in \sigma _{i}$. By definition: 
\begin{equation*}
\Phi (W_{j}^{f})=\sum_{k\in \psi _{i}}u_{jk}W_{k}^{g}\text{.}
\end{equation*} 
Since $\Phi $ is injective and $\left\{ W_{1}^{f},\ldots ,W_{n}^{f}\right\} $
and $\left\{ W_{1}^{g},\ldots ,W_{n}^{g}\right\} $ are linearly independent
sets, we have: 
\begin{eqnarray*}
\left\vert \sigma _{i}\right\vert  &=&\dim \mathrm{span}\left\{ W_{j}:j\in
\sigma _{i}\right\} =\dim \mathrm{span}\left\{ \Phi (W_{j}^{f}):j\in \sigma
_{i}\right\}  \\
&\leq &\dim \mathrm{span}\left\{ W_{k}^{g}:k\in \psi _{i}\right\}
=\left\vert \psi _{i}\right\vert \text{.}
\end{eqnarray*} 
If, moreover, $\sigma _{i}=s_{\Phi ^{-1}}\left( \psi _{i}\right) $ for all $ 
i\in \left\{ 1,\ldots ,n\right\} $, then one gets $\left\vert \sigma
_{i}\right\vert =\left\vert \psi _{i}\right\vert $ for $i\in \left\{
1,\ldots ,n\right\} $.

Now, we turn to the construction of a partition $\left\{ \sigma _{1},\ldots
,\sigma _{n}\right\} $ of $\left\{ 1,\ldots ,n\right\} $ such that $s_{\Phi
^{-1}}\circ s_{\Phi }\left( \sigma _{i}\right) =\sigma _{i}$ for all $i\in
\left\{ 1,\ldots ,n\right\} $. Let $i\in \left\{ 1,\ldots ,n\right\} $. Set $ 
A_{0}=\left\{ i\right\} $ and $A_{m+1}=s_{\Phi ^{-1}}\circ s_{\Phi }\left(
A_{m}\right) $ for all $m\in \mathbb{N}$. Let $k\in A_{m}$ for some $m\in 
\mathbb{N}$. Then there exists $j\in s_{\Phi }\left( A_{m}\right) $ such
that $u_{kj}\not=0$. Since $j\in s_{\Phi }\left( \left\{ k\right\} \right) $
we have $k\in s_{\Phi ^{-1}}\left( \left\{ j\right\} \right) \subseteq
s_{\Phi ^{-1}}\circ s_{\Phi }(A_{m})$. Hence $k\in A_{m+1}$. So $\left(
A_{m}\right) _{m\in \mathbb{N}}$ is a sequence of subsets of $\left\{
1,\ldots ,n\right\} $ increasing for the inclusion. Since $\left\{ 1,\ldots
,n\right\} $ is finite, there exists $N\in \mathbb{N}$ such that $ 
A_{N}=A_{N+1}$. We define $\mathrm{cl}\left( i\right) =A_{N}$.

We thus construct subsets $\mathrm{cl}(1),\ldots ,\mathrm{cl}(n)$ of $ 
\left\{ 1,\ldots ,n\right\} $ such that $i\in \mathrm{cl}\left( i\right) $
for all $i\in \left\{ 1,\ldots ,n\right\} $, so $\bigcup\limits_{i\in
\left\{ 1,\ldots ,n\right\} }\mathrm{cl}(i)=\left\{ 1,\ldots ,n\right\} $.
To show $\left\{ \mathrm{cl}(1),\ldots ,\mathrm{cl}(n)\right\} $ is a
partition, it is thus sufficient to show that if, for some $i,j\in \left\{
1,\ldots ,n\right\} $ $\mathrm{cl}(i)\cap \mathrm{cl}(j)\not=\emptyset $
then $\mathrm{cl}(i)=\mathrm{cl}(j)$. 

To do so, let us assume that $k\in \mathrm{cl}(i)$. Since $s_{\Phi
^{-1}}\circ s_{\Phi }\left( \mathrm{cl}\left( i\right) \right) =\mathrm{cl} 
(i)$ by definition, we conclude that $\mathrm{cl}(k)\subseteq \mathrm{cl} 
(i)$. On the other hand, by construction, there exists $j_{1},\ldots ,j_{q}$
(for $q\leq n$) such that $j_{1}=i$, $j_{q}=k$ and $j_{m+1}\in s_{\Phi
^{-1}}\circ s_{\Phi }\left( \left\{ j_{m}\right\} \right) $. Fix $m\in
\left\{ 1,\ldots ,q\right\} $. Since $j_{m+1}\in s_{\Phi ^{-1}}\circ s_{\Phi
}\left( \left\{ j_{m}\right\} \right) $ there exists $r_{m}\in s_{\Phi
}\left( \left\{ j_{m}\right\} \right) $ such that $j_{m+1}\in s_{\Phi
^{-1}}\left( \left\{ r_{m}\right\} \right) $ so $r_{m}\in s_{\Phi }\left(
\left\{ j_{m+1}\right\} \right) $. Since $r_{m}\in s_{\Phi }\left( \left\{
j_{m}\right\} \right) $ we have $j_{m}\in s_{\Phi ^{-1}}\left( \left\{
r_{m}\right\} \right) $ and therefore $j_{m}\in s_{\Phi ^{-1}}\circ s_{\Phi
}\left( j_{m+1}\right) $. Hence, $i\in \mathrm{cl}(k)$. Therefore, $ 
\mathrm{cl}(k)=\mathrm{cl}(i)$.

Assume now that $\mathrm{cl}\left( i\right) \cap \mathrm{cl} 
(j)\not=\emptyset $ for some $i,j\in \left\{ 1,\ldots ,n\right\} $. Let $ 
k\in \mathrm{cl}\left( i\right) \cap \mathrm{cl}(j)$. We have shown that $ 
\mathrm{cl}(i)=\mathrm{cl}(k)=\mathrm{cl}(j)$. Therefore, $\left\{ 
\mathrm{cl}(1),\ldots ,\mathrm{cl}(n)\right\} $ is a partition of $\left\{
1,\ldots ,n\right\} $ such that $s_{\Phi ^{-1}}\circ s_{\Phi }\left( 
\mathrm{cl}(i)\right) =\mathrm{cl}(i)$ for $i=1,\ldots ,n$. We rewrite
this partition as $\left\{ \sigma _{1},\ldots ,\sigma _{p}\right\} $ (the
order is unimportant). Setting $\psi _{i}=s_{\Phi }\left( \sigma _{i}\right) 
$ for $i=1,\ldots ,p$, we have found a partition satisfying our theorem.

Note, at last, that if $A\subseteq \left\{ 1,\ldots ,n\right\} $ such that $ 
s_{\Phi ^{-1}}\circ s_{\Phi }\left( A\right) =A$ then for $i\in A$ then 
\begin{equation*}
s_{\Phi ^{-1}}\circ s_{\Phi }(i)\subseteq s_{\Phi ^{-1}}\circ s_{\Phi
}\left( A\right) =A
\end{equation*}
and thus $\mathrm{cl}(i)\subseteq A$, as desired.
\end{proof}

We now show that the permutation-rescaling equivalence on symbols corresponds precisely to isomorphism of polynomial aspherical noncommutative domain algebras. To this end, we shall recall the following observations. Let $f=\sum_{\alpha} a_{\alpha} X_{\alpha}$ be a regular positive $n$-free formal power series.

\begin{enumerate}
\item There exists positive real numbers $b_\alpha$ for all $\alpha\in\mathbb{F}_n^+$ such that $W_j$ maps the canonical basis vector $\delta_\alpha$ at $\alpha\in \mathbf{F}_n^+$ to $\sqrt{\frac{b_\alpha}{b_{g_j \alpha}}} \delta_{g_j \alpha}$ for $j=1,\ldots,n$  \cite{Popescu08}.
\item The correspondence between the coefficients $(a_\alpha)_{\alpha\in\mathbb{F}_n^+}$ and $(b_\alpha)_{\alpha\in\mathbb{F}_n^+}$ is bijective  \cite{Popescu08}. 
\item For any $\alpha\in\mathbb{F}_n^+$ we have $\|W_\alpha^f\|=\frac{1}{b_\alpha^f}$  \cite{Popescu08}.
\item The Pythagorean identity holds for the weighted shifts $W_\alpha^f$ in the following sense:
$$\|\sum_{j=1}^k c_j W^f_{\alpha_j}\|^2 = \sum_{j-1}^n |c_j|^2 \|W^f_{\alpha_j}\|^2$$ where $c_1,\ldots,c_k$ are arbitrary complex numbers and $\alpha_1,\ldots,\alpha_n \in\mathbb{F}_n^+$ of are all words \emph{of the same length} \cite{Arias09}.
\end{enumerate}

We now can show:

\begin{theorem}\label{Main2}
Let $f,g$ be two regular, positive $n$-free polynomials. There exists an
isomorphism $\Phi :\mathcal{A}\left( \mathcal{D}_{f}\right) \longrightarrow 
\mathcal{A}\left( \mathcal{D}_{g}\right) $ such that $\widehat{\Phi }(0)=0$
if and only if $f$ and $g$ are permutation-rescaling equivalent.
\end{theorem}

\begin{proof}
Write $f=\sum a_{\alpha }^{f}X_{\alpha }$ and $g=\sum a_{\alpha
}^{g}X_{\alpha }$. Up to rescaling (see  \cite[Lemma 4.4]{Arias09}), we assume $a^f_\alpha=a^g_\alpha=1$ for all words $\alpha$ of length 1. We first assume that there exists an isomorphism $\Phi : 
\mathcal{A}\left( \mathcal{D}_{f}\right) \longrightarrow \mathcal{A}\left( 
\mathcal{D}_{g}\right) $ such that $\widehat{\Phi }(0)=0$. By Proposition
(\ref{partition}), we can find a partition $\left\{ \sigma _{1},\ldots ,\sigma _{p}\right\} 
$ of $\left\{ 1,\ldots ,n\right\} $ such that $s_{\Phi ^{-1}}\circ s_{\Phi
}\left( \sigma _{i}\right) =\sigma _{i}$ for all $i\in \left\{ 1,\ldots
,p\right\} $, where $s_{\Phi }$ (resp. $s_{\Phi ^{-1}}$) is the support
function for $\Phi $ (resp. $\Phi ^{-1}$). Up to permutation of the free
variables in $g$ we may assume that $s_{\Phi }\left( \sigma _{i}\right)
=\sigma _{i}$ for $i\in \left\{ 1,\ldots ,p\right\} $. 

Let $l\in \mathbb{N}$ with $l>1$ be given. The finite set $\left\{ b_{\alpha
}^{f},b_{\alpha }^{g}:\alpha \in \mathbb{F}_{n}^{+}\wedge \left\vert \alpha
\right\vert =l\right\} $ of real numbers has a smallest number, which we
note $b_{\omega }^{f}$ (if the minimum is reached for a coefficient for $g$
instead, we flip the notations for $f$ and $g$ for this case), with $\omega
\in \mathbb{F}_{n}^{+}$ and $\left\vert \omega \right\vert =l$. We write $ 
\omega =g_{i_{1}}\cdots g_{i_{l}}$ where $g_{1},\ldots ,g_{n}$ are the
canonical generators of $\mathbb{F}_{n}^{+}$. Now, since $\Phi $ is an
isometry and is implemented on the generators $W_{1}^{f},\ldots ,W_{n}^{f}$
by a scalar unitary $\left[ u_{ij}\right] _{1\leq i,j\leq n}$ acting on $ 
(W_{1}^{g},\ldots ,W_{n}^{g})$ (see Proposition (\ref{partition}) for notations), we have: 
\begin{eqnarray*}
\frac{1}{b_{\omega }^{f}} &=&\left\Vert W_{\omega }^{f}\right\Vert
^{2}=\left\Vert \Phi \left( W_{\omega }^{f}\right) \right\Vert ^{2} \\
&=&\left\Vert \prod\limits_{k=1}^{l}\Phi \left( W_{i_{k}}^{f}\right)
\right\Vert ^{2} \\
&=&\left\Vert \prod\limits_{k=1}^{l}\sum_{r_{k}\in s_{\Phi }\left( \left\{
i_{k}\right\} \right) }u_{i_{k}r_{k}}W_{r_{k}}^{g}\right\Vert ^{2} \\
&=&\left\Vert \sum_{r_{1}\in s_{\Phi }\left( \left\{ i_{1}\right\} \right)
}\cdots \sum_{r_{l}\in s_{\Phi }\left( \left\{ i_{l}\right\} \right)
}u_{i_{1}r_{1}}u_{i_{2}r_{2}}\cdots u_{i_{l}r_{l}}W_{r_{1}}^{g}\cdots
W_{r_{l}}^{g}\right\Vert ^{2} \\
&=&\sum_{r_{1}\in s_{\Phi }\left( \left\{ i_{1}\right\} \right) }\cdots
\sum_{r_{l}\in s_{\Phi }\left( \left\{ i_{l}\right\} \right) }\left\vert
u_{i_{1}r_{1}}\right\vert ^{2}\left\vert u_{i_{2}r_{2}}\right\vert
^{2}\cdots \left\vert u_{i_{l}r_{l}}\right\vert ^{2}\frac{1}{ 
b_{g_{r_{1}}\cdots g_{r_{l}}}^{g}}\text{.}
\end{eqnarray*} 
Thus, since $U$ is unitary, we conclude that $\frac{1}{b_{\omega }^{f}}$ is
a convex combination of elements in $\left\{ \frac{1}{b_{\alpha }^{f}},\frac{1}{b_{\alpha
}^{g}}:\alpha \in \mathbb{F}_{n}^{+}\wedge \left\vert \alpha \right\vert
=l\right\} $, though it is its maximum. This is only possible if: 
\begin{equation*}
\forall r_{1}\in s_{\Phi }\left( \left\{ i_{1}\right\} \right) \ \ \cdots \
\ \forall r_{l}\in s_{\Phi }\left( \left\{ i_{1}\right\} \right) ~~~~~\frac{1 
}{b_{g_{r_{1}}\cdots g_{r_{l}}}^{g}}=\frac{1}{b_{\omega }^{f}}\text{.}
\end{equation*}

Let us adopt the following notation. A word $\alpha $ will be of type $ 
\sigma _{i_{1}}\cdots \sigma _{i_{l}}$ if $\alpha =g_{r_{1}}\cdots g_{r_{n}}$
for $r_{j}\in \sigma _{j}$ with $j=1,\ldots ,l$. Let $\sigma _{v_{1}}\cdots
\sigma _{v_{l}}$ be the type of the $\omega $ as above. By repeating the
above argument, we can thus show that $b_{\alpha }^{f}=b_{\alpha
}^{g}=b_{\omega }^{f}$ for all $\alpha $ of type $\omega $.

We can now repeat the proof above by picking the minimum of 
\begin{equation*}
\left\{ \frac{1}{b_{\alpha }^{f}},\frac{1}{b_{\alpha }^{g}}:\alpha \in 
\mathbb{F}_{n}^{+}\wedge \left\vert \alpha \right\vert =l\wedge \alpha \text{
is not of the same type as }\omega \right\} 
\end{equation*} 
and so forth until we exhaust the set $\left\{ \frac{1}{b_{\alpha }^{f}}, 
\frac{1}{b_{\alpha }^{g}}:\alpha \in \mathbb{F}_{n}^{+}\wedge \left\vert
\alpha \right\vert =l\right\} $ to show that $b_{\alpha }^{f}=b_{\alpha
}^{g}=b_{\omega }^{f}$ for all $\alpha \in \mathbb{F}_{n}^{+}$ with $ 
\left\vert \alpha \right\vert =l$. This completes our proof.
\end{proof}

\bigskip 

\begin{corollary}
Let Let $f,g$ be two regular positive $n$-free polynomials, and assume $f$
is aspherical. Then $\mathcal{A}\left( \mathcal{D}_{f}\right) $ and $ 
\mathcal{A}\left( \mathcal{D}_{g}\right) $ are isomorphic if and only if $f$
and $g$ are scale-permutation equivalent.
\end{corollary}

\begin{proof}
This results follows from Theorem (\ref{main}) and Theorem (\ref{Main2}).
\end{proof}

\bigskip 

We can summarize the current understanding of classification for
noncommutative domain algebra:

\begin{theorem}
Let $f$ be a regular, positive $n$-free polynomial. Then:

\begin{itemize}
\item If $f=\sum c_{i}X_{i}$ for some $c_{1},\ldots ,c_{n}\in \left(
0,\infty \right) $ then $\mathcal{A}\left( \mathcal{D}_{f}\right) $ is
isomorphic to the disk algebra $\mathcal{A}_{n}$,

\item If $f$ is aspherical, then $\mathcal{A}\left( \mathcal{D}_{g}\right) $
is isomorphic to $\mathcal{A}\left( \mathcal{D}_{f}\right) $ if and only if $ 
f=g$ after rescaling/permutation of the free variables of $g$.
\end{itemize}
\end{theorem}

This leaves the matter of classifying $\mathcal{A}\left( \mathcal{D} 
_{f}\right) $ when $f$ is spherical, i.e. $\mathcal{D}_{f}\left( \mathbb{C} 
\right) $ is the closed unit ball of $\mathbb{C}^{n},$ yet $f$ is not of
degree 1.

\section{Cartan's Lemma}

In \cite[Theorem 1.4]{Popescu10} and \cite[Theorem 4.5]{Popescu11},
Popescu establishes a generalization of Cartan's Lemma, first in the context
of the unit ball of $\mathcal{B}(\mathcal{H})$ for a Hilbert space $\mathcal{H}$, then for a large class of noncommutative domains that includes the domains
considered in this paper. We propose to illustrate in this section that our
methods, as implemented in this paper and in \cite{Arias09}, can be used to
obtain a simpler proof of these results. We refer to \cite{Popescu08,Popescu10,Popescu11} for definitions and the general theory of holomorphic functions in the context of noncommutative domains. We shall only use the following special case of holomorphic functions:

\begin{definition}
Let $f$ be a positive, regular $n$-free formal power series. A holomorphic map $F$ with domain and codomain $\mathcal{D}_f$ is the given of a family of complex coefficients $(c_\alpha)_{\alpha\in\mathbb{F}^+_n}$ such that, for any Hilbert space $\mathcal{H}$, the function:
$$F_{\mathcal{H}} :  (T_1,\ldots,T_n)\in\mathrm{interior}\, \mathcal{D}_f(\mathcal{H}) \longmapsto \left(\sum_{\alpha \in \mathbb{F}_n^+} c_{1,\alpha} T_\alpha,\ldots,\sum_{\alpha \in \mathbb{F}_n^+} c_{1,\alpha} T_\alpha\right)$$ is well-defined and its range is a subset of the interior of $\mathcal{D}_f(\mathcal{H})$.
\end{definition}

\begin{remark}
The domain and codomain of a holomorphic function, in this context, is the noncommutative domain $\mathcal{D}_f$ which, itself, is not a set, but a map from separable Hilbert spaces to subsets of the algebra of linear bounded operators on the given Hilbert space. Note moreover that the domain and codomain of the maps induced on various noncommutative domains by holomorphic maps are the \emph{interior} of the noncommutative domains.
\end{remark}

\begin{definition}
Let $f$ be a positive, regular, $n$-free formal power series. A biholomorphic map $F$ on $\mathcal{D}_f$ is a holomorphic map from $\mathcal{D}_f$ to $\mathcal{D}_f$ such that there exists a holomorphic map $G$ from $\mathcal{D}_f$ to $\mathcal{D}_f$ such that $F\circ G$ and $G\circ F$ both induce the identity on $\mathcal{D}_F(\mathcal{H})$ for all Hilbert spaces $\mathcal{H}$.
\end{definition}

\begin{remark}
By definition, the map induced by a holomorphic map on a specific noncommutative domain is a holomorphic map on the interior of this domain.
\end{remark}

\begin{theorem}[Cartan's Lemma]
Let $f$ be a regular positive $n$-free formal power series. Let $F$ be a biholomorphic map of ${D}_f$ such that $F(0)=0$. Then $F$ is linear.
\end{theorem}

\begin{proof}
By definition, there exists complex coefficients $(c_{j,\alpha})_{j\in\{1,\ldots n\},\alpha\in\mathbb{F}_n^+}$ such that for all Hilbert space $\mathcal{H}$ and all $(T_1,\ldots,T_n)\in\mathrm{interior}\, \mathcal{D}_f(\mathcal{H})$ we have:
$$F_{\mathcal{H}}(T_1,\ldots,T_n) = \left(\sum_{\alpha \in \mathbb{F}_n^+} c_{1,\alpha} T_\alpha,\ldots,\sum_{\alpha \in \mathbb{F}_n^+} c_{1,\alpha} T_\alpha\right)\in\mathrm{interior}\,\mathcal{D}_f(\mathcal{H})\textrm{.}$$

When $\mathcal{H}=\mathbb{C}^k$, we will write $F_k$ for the biholomorphic map induced by $F$ on $\mathbb{D}_f^k$ for $k\in\mathbb{N}$, $k>0$.

\bigskip Since this proof is essentially the same as \cite[Theorem 3.18]{Arias09}, we shall only deal with the case where $\alpha\in\mathbb{F}_n^+,|\alpha|\leq 2$ and $n=2$. 

By definition, $F$ induces a biholomorphic map $F_1$ on $\mathbb{D}_f^1$. This map fixes $0$ and  $\mathbb{D}_f^1$ is a circular domain (even Reinhardt) in $\mathbb{C}^n$, Cartan's lemma \cite{Krantz} implies that $F_1$ is linear.

Now, the $j^{th}$ coordinate of $F_1(z_1,z_2)$ for an arbitrary $(z_1,z_2)\in\mathbb{D}_f^1$ is:
$$c_{j,g_1} z_1 + c_{j,g_2} z_2 + c_{j,g_1 g_1} z_1 ^2 + c_{j,g_2 g_2} z_2 ^2 + (c_{j,g_1 g_2} + c_{j,g_2 g_1}) z_1 z_2 + \cdots $$
so the linearity of $F_1$ implies that:
$$ c_{j,g_1 g_1} z_1 ^2 + c_{j,g_2 g_2} z_2 ^2 + (c_{j,g_1 g_2} + c_{j,g_2 g_1}) z_1 z_2 + \cdots = 0$$
which in turns shows that $c_{j,g_1 g_1} = c_{j,g_2 g_2} = c_{j,g_1 g_2} + c_{j,g_2 g_1} = 0$.
To show that $c_{j,g_1 g_2}=c_{j,g_2 g_1} = 0$ we go to higher dimensions. Again by definition, $F$ induces a biholomorphic map on $\mathbb{D}_f^2$. The later domain is circular (not Reinhardt in general), and thus again by Cartan's lemma \cite{Krantz}, since $F_2(0)=0$ we conclude that $F_2$ is linear.

Now, a quick computation shows that the $(2,1)$ component of the $j^{th}$ coordinate of $F_2(M,N)$ where $M=\left[ \begin{array}{cc} \lambda_1 & \lambda_2 \\ \lambda_3 & \lambda_4 \end{array} \right]$ and $N=\left[ \begin{array}{cc} \lambda_5 & \lambda_6 \\ \lambda_7 & \lambda_8 \end{array} \right]$  are $2\times 2$ complex matrices in $\mathbb{D}_f^2$ is given by:

\begin{eqnarray*}
 c_{j,g_1} \lambda_2 + c_{j,g_2} \lambda_6 &+& c_{j,g_1 g_2} (\lambda_1\lambda_6 + \lambda_2 \lambda_8) + c_{j,g_2 g_1} (\lambda_2 \lambda_5 + \lambda_4 \lambda_6 )\\
&+& \textrm{terms of higher degrees in }\lambda_1,\ldots,\lambda_8
\end{eqnarray*}
so linearity of $F_2$ implies that $c_{j,g_1 g_2}=c_{j,g_2 g_1} =0$.

The proof for higher terms is similar and undertaken in \cite[Theorem 3.18]{Arias09}.

\end{proof}

We expect that the same method of reduction to finite dimension can yield other generalizations of results from the study of domains in complex analysis to the framework of Popescu's noncommutative domains.

\bibliographystyle{amsplain}
\bibliography{../thesis}

\end{document}